\documentclass{article}
\usepackage[]{epsf,epsfig,amsmath,amssymb,amsfonts,latexsym}
\title{Extremal points of  high dimensional random walks and mixing times of a brownian motion on the sphere}
\author{Ronen Eldan\thanks{Partially supported by the Israel Science Foundation and by a Farajun Foundation Fellowship.}}
\newtheorem{theorem}{Theorem}[section]
\newtheorem{question}{Question}[section]
\newtheorem{lemma}{Lemma}[section]
\newtheorem{proposition}{Proposition}[section]

\newtheorem{corollary}{Corollary}[section]

\newtheorem{remark}{Remark}[section]

\def \PP {\mathbb P}
\def \P {\mathbb P}
\def\qed{\hfill $\vcenter{\hrule height .3mm
\hbox {\vrule width .3mm height 2.1mm \kern 2mm \vrule width .3mm
height 2.1mm} \hrule height .3mm}$ \bigskip}
\def\P{\mathbb{P}}

\def\EE{\mathbb{E}}

\def\RR{\mathbb{R}}

\def\Sph{S^{n-1}}

\begin{document}
\maketitle

\begin{abstract}
We derive asymptotics for the probability that the origin is an extremal point of a random walk in $\RR^n$. We show that in order for the probability to be roughly $1/2$, the number of steps of the random walk should be between $e^{n/(C \log n)}$ and $e^{C n \log n}$ for some constant $C>0$. As a result, we attain a bound for the $\frac{\pi}{2}$-covering time of a spherical brownian motion.
\end{abstract}
\bigskip
\section{Introduction}

The object of this paper is to address the following question: given a random walk in Euclidean space, how long does it typically take until the starting point of the random walk ceases to be an extremal point of its range? We approach this question from a high-dimensional point of view. In particular, we try to derive asymptotics of some quantities related to this question, as the dimension goes to infinity. \\

Let us give a more precise formulation of our question. Fix a dimension $n \in \mathbb{N}$. For a set $K \subset \RR^n$, by $\partial K$ we denote its  boundary, by $Int(K)$ its interior, and by $conv(K)$ we denote its convex hull. Let $t_1 \leq ... \leq t_N$ be a Poisson point process on $[0,1]$ with intensity $\alpha$, let $B(t)$ be an $n$-dimensional standard brownian motion. Define $X_0 = 0, X_i = B(t_i)$. We call $X_1,...,X_N$ a random walk in $\RR^n$. We say that the origin is an extremal point of this random walk if $0 \in \partial K$, where $K := conv(\{X_0, X_1,...,X_N\})$. \\ \\
Denote by $p(n,\alpha)$ the probability that the origin is an extremal point of the
the random walk $X_0, X_1,...,X_N$. For $n \in \mathbb{N}$, note that $p(n,\alpha)$ is a decreasing function of $\alpha$ and denote by $\alpha(n)$ the smallest number, $\alpha$, such that $p(n, \alpha) \leq \frac{1}{2}$. Our aim in this note is to prove the following asymptotic bound:
\begin{theorem} \label{mainttt}
With $\alpha(n)$ defined as above, one has
$$
e^{c n / \log n} < \alpha(n) < e^{C n \log n}.
$$
for some universal constants $c,C>0$.
\end{theorem}

\medskip
Following rather similar lines, one can also prove that the same asymptotics
are correct for the standard random walk on $\mathbb{Z}^n$. Namely, one can prove the following result:
\begin{theorem} \label{discrete}
Let $S_1,...,S_N$ be the standard random walk on $\mathbb{Z}^n$. Define,
$$N(n) = \min \left \{ N \in \mathbb{N} ~~ \left | ~~ \P \left (0 \mbox{ is an extremal point of } conv \{S_1,...,S_N\} \right ) \leq \frac{1}{2} \right .  \right \}$$
Then,
$$
e^{c n / \log n} < N(n) < e^{C n \log n}.
$$
for some universal constants $c,C>0$.
\end{theorem}

The latter theorem may be, in fact, more interesting for probabilists than the former. Nevertheless, we choose to omit some of the details of its proof since it is more involved than the proof of theorem \ref{mainttt}, and the two proofs share the same ideas. We will provide an outline of proof along with some remarks about the further technical work that should be done in order to prove it. 

\begin{remark}
By means of the so-called \emph{reflection principle}, it may be shown that for a $1$-dimensional, simple random walk,
the probability to remain non-negative after $N$ steps is of the order $1 / \sqrt{N}$. The expectation of the first
time it becomes negative is therefore $\infty$. It follows that the expectation of the first time that the convex hull
of a random walk in any dimension contains the origin in its interior is also infinite.
\end{remark}
\bigskip

A corollary of the above result concerns with covering times of the spherical brownian motion. We define $\Sph = \{x \in \RR^n, |x|=1 \}$, $| \cdot |$ being the standard Euclidean norm. Given a standard brownian motion $B(t)$ in $\RR^n$, $n>2$, the function $\theta(t) = \frac{B(t)}{|B(t)|}$ is almost surely defined for all $t>0$. By the Dambis / Dubins-Schwarz theorem, there exists a non-decreasing (random) function $T(\cdot)$ such that $\theta(T(\cdot))$ is a strong Markov process whose quadratic variation as time $t$ is equal to $(n-1)t$. We refer to the process $\theta(T(t))$ as a spherical brownian motion (or a brownian motion on $\Sph$). Furthermore, we denote by $d(\cdot, \cdot)$ the geodesic distance on $\Sph$, equipped with the standard metric. The $\epsilon$-neighbourhood of a point $x \in \Sph$ is defined as $\nu_x(\epsilon) = \{y \in \Sph, ~d(x,y) < \epsilon \}$. We say that a set $A \subset \Sph$ is an $\epsilon$-covering of the sphere if $\bigcup_{x \in A} \nu_x(\epsilon) = \Sph$. \\

Let us now consider the following question: given a brownian motion on $\Sph$, how long does it typically take until the path is not contained in an open hemisphere? Equivalently, how long does it take for a brownian motion to be a $\pi / 2$-covering of the sphere? Covering times of random walks and brownian in different settings is a subject that has been widely studied in the past decades (see e.g., \cite{Ad}, \cite{DPRZ}, \cite{M} and references therein). Matthews \cite{M} studied the $\epsilon$-cover time for brownian motion on an $n$-dimensional sphere. In his work, he considers the asymptotics as $\epsilon$ tends to zero and the dimension is fixed. \\

One motivation for the study of covering times on the sphere is a technique for viewing multidimensional data developed by Asimov \cite{A}, known as the Grand Tour. In this technique, a high dimensional object (usually, a measure on $\RR^n$) is analyzed through visual inspection of its projections onto subspaces of small dimension. When considering one-dimensional marginals, the set of directions may be taken from the range of a spherical brownian motion. In this case, one may be interested in estimating how long should takes for the brownian motion to visit the a certain neighbourhood of all possible directions on the sphere, thus indicating that the set of inspected marginals is rather dense. \\

Let $E(n)$ be the expected time it takes the spherical brownian is a $\frac{\pi}{2}$-covering of the sphere, in other words,
$$
E(n) = \EE \left [\inf \left \{t > 0 ;~ 0 \mbox{ is in the interior of } conv(\{\mathbf{SB}_n(s); ~ 0 < s \leq t \}) \right  \} \right ],
$$
where $\mathbf{SB}_n(s)$ is brownian motion on $\Sph$. A corollary of our bounds for $\alpha(n)$ is a corresponding bound for the asymptotics of $E(n)$, as $n$ goes to infinity. Namely,

\begin{corollary} \label{covering}
There exists a universal constant $C>0$ such that,
$$ 
\frac{1}{C \log n} < E(n) < C \log n, ~~~ \forall n \geq 1.
$$
\end{corollary}

The above corollary and the work of Matthews complete each other in a certain sense: The asymptotics derived by Matthews for $E(n)$ in the case of $\epsilon$-covering, when $\epsilon \to 0$, is roughly $E(n) \sim \sqrt{n} \epsilon^{n-3} \log(\epsilon^{-1})$. In other words, for small $\epsilon$, the
time is exponential in the dimension. Our result therefore suggests a rather significant phase shift as $\epsilon$ approaches $\pi / 2$.

Another possible application the last corollary is related to the following illumination problem: a high dimensional convex object (say, a planet) is rotating randomly. A single light source is located very far from the object. How long will it take until every point on the surface of the object has been illuminated at least once? 

\medskip
The organization of the rest of this paper is as follows: the lower bound of theorem \ref{mainttt} will be proven in section \ref{secoelower} and the upper bound will be proven in section \ref{secoeupper}. Section \ref{secoedisc} is devoted to filling some of the missing details for the proof of theorem \ref{discrete}. In section 5, we prove corollary \ref{covering}. Finally, in section \ref{secoeremarks}, we list some further facts that can be derived using the same methods of proof and raise some questions for possible further research.  \\

Throughout this note, the symbols $C,C',C'',c,c',c''$ denote positive universal constants whose values may change between different formulas. We write
$f(n) = O(g(n))$ if there is a positive constant $M>0$ such that $f(n) < M(g(n))$ for all $n$, and we write $f(n) = o(g(n))$ if $f(n) / g(n) \to 0$ as $n \to \infty$. Given a subset $A \subset \RR^n$, by $conv(A)$ we denote the convex hull of $A$. Given two random variables, $X$ and $Y$, the notation $X \sim Y$ is to say that the two variables have the same distribution. For random vector $X \in \RR^n$ we denote its barycenter by $b(X) := \EE[X]$, and its covariance matrix by
$Cov(X) := \EE[(X - b(X)) \otimes (X - b(X))]$. \\ \\

\textbf{Acknowledgements} I would like to express my thanks to Itai Benjamini for introducing me to the question and for several discussions about the subject, and to Bo`az Klartag, Boris Tsirelson and Ron Peled for useful discussions. Finally, I would like to thank the referee of this paper for numerous enlightening comments and useful suggestions.

\section{The Lower Bound} \label{secoelower}
The aim of this section is to prove the following bound:
\begin{theorem} \label{lower}
There exists a universal constant $c>0$ such that the following holds: Suppose $\alpha < e^{c n / \log n}$. 
Let $B(t)$ be a standard brownian motion in $\RR^n$. Then,
\begin{equation} \label{eqlower1}
\P \left ( 0 \mbox{ is in the interior of } conv( \{B(t) ~ | ~ \alpha^{-1} \leq t \leq 1 \}) \right ) < 0.1.
\end{equation}
In particular, if $t_1 \leq ... \leq t_N$ are points generated according to a poisson process on $[0,1]$ with intensity $c \alpha$, independently of
$B(t)$, then
\begin{equation} \label{eqlower2}
\P \left ( 0 \mbox{ is an extremal point of the set } \{B(0), B(t_1),...,B(t_N) \} \right ) > \frac{1}{2}.
\end{equation}
\end{theorem}
\bigskip
Before we begin the proof, we will need the following ingredient: recall Bernstein's inequality, \cite{U}, which can be states as follows.
\begin{theorem} \label{berns} (Bernstein's inequality) 
Let $X_1,...,X_n$ be independent random variables. Suppose that for
some positive $L>1$ and every integer $k>0$,
\begin{equation} \label{bernsineq}
\EE[|X_i - \EE[X_i]|^k] < \frac{\EE[X_i^2]}{2} L^{k-2} k!
\end{equation}
Then,
$$
\P \left ( \left | \sum_{i=1}^n (X_i - \EE[X_i]) \right | > 2t \sqrt {\sum_{i=1}^n Var[X_i]} \right ) <e^{-t^2}
$$
for every $0<t<\frac{\sqrt{\sum_{i=1}^n Var[X_i]}}{2L}$.
\end{theorem} 
\medskip
\textbf{Proof of theorem \ref{lower}:} \\

First of all, we note that equation (\ref{eqlower2}) follows easily from equation (\ref{eqlower1}). Indeed, by a small enough choice of the constant $c$, we can make sure that with probability at least $3/4$, none of the points $t_1,...,t_N$ fall inside the interval $[0,\alpha^{-1}]$. We turn 
to prove equation (\ref{eqlower1}). \\ \\
By choosing a suitable (small enough) value for the constant $c$, we may always assume that the dimension, $n$, is larger than some universal constant.
Define $m = \left \lfloor \frac{c n}{\log n} \right \rfloor$, where the value of the constant $c>0$ will be chosen later on. Since the probability in equation (\ref{eqlower1}) is increasing with $\alpha$, we may assume that $\alpha=2^{m-1}$. Moreover, in order to simplify the below formulas, we note that by using a scaling argument we can assume that our time interval is $[0,2^{m-1}]$ (rather than the interval $[0,1]$), and show that,
$$
\P \left ( 0 \mbox{ is in the interior of } conv(\{B(t) ~ | ~ 1 \leq t \leq 2^{m-1} \}) \right ) < \frac{1}{4}.
$$
We will show that with high probability there exists a vector $v$ which demonstrates that the origin is not in the interior, i.e, that $\langle B(t), v \rangle > 0$ for all $1 \leq t \leq 2^{m-1}$. \\ \\
The construction of the vector $v$ is as follows. Define,
$$
v_i = B \left (2^{i} \right )- B \left (2^{i-1} \right ),
$$
for $i=0,...,m-1$, and
$$
v = \frac{1}{\sqrt m} \sum_{i=0}^{m-1} \frac{v_i}{\sqrt{\EE[|v_i|^2]}} =  \frac{1}{\sqrt m} \sum_{i=0}^{m-1} \frac{v_i}{ \sqrt n (\sqrt{2})^{i-1}}.
$$ 
Note that the vectors $\frac{v_i}{\sqrt{\EE[|v_i|^2]}}$ are independent, identically distributed gaussian random vectors with expectation $0$
and with covariance matrix $\frac{1}{n} Id$. It follows that the vector $v$ is also a gaussian random vector whose expectation is $0$ and whose covariance matrix is equal to $\frac{1}{n} Id$. A calculation then gives,
\begin{equation} \label{normv}
\PP \left (\frac{1}{2} < |v| < 2 \right ) > 1 - e^{-c' n}
\end{equation}
for some universal constant $c'>0$. \\ \\
Fix $0 \leq k \leq m-1$. Let us inspect the scalar product $p = \langle B(2^k), v \rangle$. 
for all $0 \leq i \leq m-1$, we denote $v_i = (v_{i,1},...,v_{i,n})$. Note that both $B(2^k)$ and $v$ are linear combinations of $v_i$'s with deterministic coefficients, hence $p$ admits the form
$$
p = \sum_{j=1}^n \sum_{i=0}^{m-1} \sum_{l=0}^{m-1} \alpha_i \beta_l v_{i,j} v_{l,j}
$$
for some constants $\{ \alpha_i \}_{i=0}^{m-1}, \{ \beta_l\}_{l=0}^{m-1}$. Define,
$$
w_j = \sum_{i=1}^m \sum_{l=1}^m \alpha_i \beta_l v_{i,j} v_{l,j}, ~~ \mbox{for } j=1,..,n.
$$
Clearly, the $w_j$'s are independent and identically distributed, so there exist numbers $a,b$ such that
\begin{equation} \label{wrep}
w_j \sim X (a X + bY)
\end{equation}
where $X,Y$ are independent standard gaussian random variables.  \\
Our next goal is to calculate the expectation and the variance of $w_j$. To that end, we may write, for all $j=1,..,n$,
\begin{equation} \label{uglysum}
w_j = \left ( \sum_{i=0}^{k} v_{i,j} \right ) \left (\frac{1}{\sqrt {n m}} \sum_{l=0}^{m-1} \frac{v_{l,j}}{(\sqrt{2})^{l-1}}  \right ) = 
\frac{1}{\sqrt {n m}} \sum_{i=0}^k \sum_{l=0}^{m-1} \frac{1}{(\sqrt{2})^{l-1}} v_{i,j} v_{l,j}.
\end{equation}
So,
$$
\EE[w_j] \geq \frac{1}{\sqrt{nm}} \frac{\EE[v_{k,j}^2]}{\sqrt{2}^{k-1}} = \frac{\sqrt{2}^{k-1} }{\sqrt{nm}},
$$
which means that,
\begin{equation} \label{expest}
\EE[p] \geq \frac{\sqrt{2}^{k-1} \sqrt{n}}{\sqrt{m}}.
\end{equation}
Next, in order to estimate $Var[w_j]$ we use (\ref{uglysum}) again to obtain,
$$
\EE[w_j^2] = \frac{1}{nm} \EE \left [ \left ( \sum_{i=0}^k \sum_{l=0}^{m-1} \frac{1}{(\sqrt{2})^{l-1}} v_{i,j} v_{l,j} \right )^2 \right ] =
$$
$$
\frac{1}{nm} \left ( \sum_{{i \neq l, 0 \leq i \leq k}, \atop {0 \leq l \leq m-1}} \frac{1}{2^{l-1}} \EE[ v_{l,j}^2] \EE[ v_{i,j}^2] + \sum_{{i \neq l,} \atop {0 \leq i,l \leq k}} \frac{1}{\sqrt{2}^{i+l-2}} \EE[v_{l,j}^2]\EE[v_{i,j}^2] + \sum_{i=0}^k \frac{1}{2^{i-1}} \EE [v_{i,j}^4] \right ) \leq
$$
$$
\frac{1}{nm} \left ( m \sum_{i=0}^k 2^i  + 2 \sum_{0 \leq i \leq l \leq k} \frac{1}{2^{i-1}} \EE[v_{l,j}^2]\EE[v_{i,j}^2] + 3 \sum_{i=0}^k 2^i \right ) < \frac{2^{k+2}}{n}.
$$
So,
\begin{equation} \label{varest}
Var[p] < 2^{k+2}.
\end{equation}
Note that $\EE[p] > \sqrt {\frac{n}{8m} Var[p]} > \sqrt{0.1 c^{-1} \log n}\sqrt{ Var[p]} $. \\ \\
It follows from representation (\ref{wrep}), from that fact that a standard Gaussian random variable, $X$, satisfies
$\EE[|X|^p] \leq p^{p/2}$ for all $p>1$, and from the Cauchy-Schwartz inequality that,
\begin{equation} \label{moments}
\EE[|w_j - \EE[w_j]|^p] < (10 Var[w_j])^{p/2} p!, ~~ \forall p \in \mathbb{N}.
\end{equation}
We may therefore invoke theorem \ref{berns} on the random variables $w_j$. Setting $t=\sqrt{\frac{n}{10m}}$, $L=10 \frac{\sqrt{2}^{k+2}}{\sqrt n}$ and plugging into (\ref{bernsineq}) leads to:
$$
\P \left (|p - \EE[p]| > \sqrt{\frac{m}{10 n}} \sqrt{Var[p]} \right ) < e^{-\frac{n}{10 m}}.
$$
Plugging in (\ref{expest}) and (\ref{varest}) and using the assumption that $c$ can be smaller than any universal constant gives, 
$$
\P(p < \frac{1}{2} \EE[p]) < e^{-\frac{n}{10 m}} < n^{-5}.
$$
Define $A$ to be the following event:
$$
A = \left \{  \langle v, B(2^k) \rangle > \frac{1}{2} \sqrt{\frac{n}{m}} \sqrt{2}^{k-1}, ~~ \forall 0 \leq k \leq m - 1  \right \}
$$
Applying a union bound for $k=0,...,m-1$, we learn that 
\begin{equation} \label{pofa}
\P(A) > 1 - \frac{1}{n^2}.
\end{equation}
Recall that the distribution of the maximal value of a brownian bridge (see e.g., \cite{sw}, page 34) starting at $y=a$ at time $0$ and ending at $y=b$ at time $T$ is,
\begin{equation} \label{bridgeprop}
f_{M^{a \to b}(T)}(y) = \mathbf{1}_{\{ y \notin [a,b] \}} 4 \frac{y - \frac{a+b}{2}}{T} e^{-\frac{2}{T}(y-a)(y-b)}.
\end{equation}
Define the events,
$$
C_k := \{\langle B(t), v \rangle > 0, ~~ \forall 2^k \leq t \leq 2^{k+1}\}.
$$
Our next goal is to show that when conditioning on $A$, the probability of $C_k$ is close to one, using the following idea:
instead of generating the brownian motion, one can alternatively generate the points $B(2^k)$ and then "fill in" the missing gaps by independent brownian bridges. When the event $A$ holds, the endpoints of the bridges $\langle B(t), v \rangle, ~ 2^k \leq t \leq 2^{k+1}$ are quite large with
respect to the standard deviation of their midpoint, and we may use (\ref{bridgeprop}). \\
More formally, Let $\tilde B(t)$ be a brownian bridge such that $B(0)=B(1)=0$, independent of $B(t)$. Define, 
$$
B_k(t) = B(2^k) + (B(2^{k+1}) - B(2^{k})) t + \sqrt{2}^k \tilde B(t).
$$ 
By a representation theorem for the brownian bridge, the functions $B_k(t)$ and $B(2^k + 2^k t)$ share the same distribution. Moreover, if an event $\tilde A$ is measurable by the sigma algebra generated by the points $B(2^j), 0 \leq j \leq m-1$, then the distribution of these two functions
is the same, event when conditioned on the event $\tilde A$. Therefore, one has,
$$
\P(C_k | A) = \P(\langle B_k(t), v \rangle > 0, ~~ \forall 0 \leq t \leq 1 ~ | ~ A).
$$
Since the maximum of a brownian bridge is monotone with respect to its endpoints, it follows that
\begin{equation} \label{bridges}
\P(\langle B_k(t), v \rangle > 0, ~~ \forall 0 \leq t \leq 1 | A) > \P \left (\langle \tilde B(t), v \rangle < \sqrt{\frac {n} {8m}}, ~~ \forall 0 \leq t \leq 1 \right ).
\end{equation}
Using (\ref{bridgeprop}) then yields,
\begin{equation} \label{bridgebound}
\P(C_k ~|~A) > 1 - \exp \left (- \log n / (8 c |v|^2) \right ).
\end{equation}
Using the above with (\ref{normv}) and choosing $c$ small enough, we get
$$
\P(C_k ~|~A) > 1 - \frac{1}{n^3}.
$$
Finally, combining with (\ref{pofa}) and using a union bound yields,
$$
\P \left (\langle B(t), v \rangle > 0, ~~ \forall 1 \leq t \leq 2^{m-1} \right ) > \P(A) \left (1- \sum_{k=1}^m (1 - \P(C_k | A)) \right ) > 1 - \frac{1}{n}.
$$
The proof is complete. \qed

\section{The Upper Bound} \label{secoeupper}

The goal of this section is the proof of the following estimate:
\begin{theorem} \label{upper}
There exists a universal constant $C>0$ such that the following holds: Let $\alpha = e^{C n \log n}$.  Let $t_1 \leq ... \leq t_N$  be points generated according to a poisson process on $[0,1]$ with intensity $\alpha$, and let $B(t)$ be a standard brownian motion, independent of the point process. Consider the random walk $B(0), B(t_1),...,B(t_N)$.  The probability that the origin is an extremal point of this random walk is smaller than $n^{-n}$.
\end{theorem}

We open the section with some well-known facts concerning the probabilities that random
walks and discrete brownian bridges stay positive. Again let $0 \leq t_1 \leq ... \leq t_N \leq 1$ be a poisson 
point process on $[0,1]$ with intensity $\alpha$, and let $W(t)$ be a standard 1-dimensional
brownian motion. Consider the random walk  $W(0), W(t_1),...,W(t_N)$. By slight abuse of notation, for $1 \leq j \leq n$, denote $W(j) = W(t_j)$. 
Let us calculate the probability that $W(j) \geq 0$ for all $1 \leq j \leq N$. \\
Recall the second arcsine law of P.Levi, (see for example \cite{MP}, Chapter 5, p. 137). Define a random variable,
$$
X = \int_0^1 \mathbf{1}_{\{ W(t) < 0 \}} dt.
$$
According to the second arcsine law, $X$ has the same distribution as $(1 + C^2)^{-1}$ where $C$ is a Cauchy 
random variable with parameter 1. Using the definition of the Poisson distribution, this means that,
$$
\P(B(t_i) > 0, ~~ \forall 1 \leq i \leq N(m)) = \EE \left [e^{-\alpha (1 + C^2)^{-1}} \right ] = \frac{1}{\pi} \int_{- \infty}^{\infty} e^{-\frac{\alpha}{1 + x^2}} \frac{1}{1 + x^2} dx = 
$$
$$
\frac{2}{\pi} \int_0^{\pi/2} e^{-\alpha \cos^2 t} dt = \frac{1}{\pi} \int_0^1 e^{-\alpha w} \frac{1}{\sqrt{w(1-w)}} dw = 
$$
$$
\frac{1}{\pi \sqrt{\alpha}} \int_0^\alpha e^{-s} \frac{1}{\sqrt{s(1 - \frac{s}{\alpha})}} ds.
$$
It is easy to check that the latter integral has a limit as $\alpha \to \infty$. Consequently,
\begin{equation} \label{walkest}
\P(B(t_i) > 0, ~~ \forall 1 \leq i \leq N) = \frac{1}{\sqrt \alpha} \left (\frac{1}{\pi} \int_0^\infty \frac{e^{-s}}{\sqrt s} ds \right ) \left (1 + o \left (\frac{1}{\alpha} \right ) \right ) = 
\end{equation}
$$
\frac{1}{\sqrt {\pi \alpha}} \left (1 + o \left (\frac{1}{\alpha} \right ) \right ).
$$
$$~$$
Now suppose that $W(t)$ is a brownian bridge such that $W(0)=W(1)=0$ and consider the discrete brownian bridge 
$W(0), W(t_1),...,W(t_N), W(1)$. 

The cyclic shifting principle (see e.g., \cite{B}) is the following observation: for every $0 \leq s \leq 1$, define
$\Gamma_s(t) = t + s$, where the sum is to be understood as a sum on the torus $[0,1]$. Then the function $W \circ \Gamma_s (t) - W(s)$ has 
the same distribution as the function $W(t)$. Now, since there is exactly one choice $i$ between $1$ and $N$ such that
$W(t_j) - W(t_i)$ will be non-negative for every $1 \leq j \leq N$, it follows that for only one choice of $1 \leq i \leq N$, the function
$$
W \circ \Gamma_{t_i}(\cdot) - W(t_i)
$$
will be positive for all the points $t_j - t_i$, $1 \leq j \leq N$ (where the subtraction is again understood on the torus $[0,1]$). Since
the points $t_1,...,t_N$ are independent of the function $W(t)$, it follows that
\begin{equation} \label{bridgeest}
\PP(W(t_i) \geq 0, ~~ \forall 1 \leq i \leq N) = \EE \left [\frac{1}{N} \right ] = \frac{1}{\alpha} + O \left (\frac{1}{\alpha^{3/2}} \right ).
\end{equation}
(recall that $N$ was a poisson random variable with expectation $\alpha$).  \\ \\
We now have the necessary ingredients for proving the upper bound. \\ \\
\textbf{Proof of theorem (\ref{upper})}: \\
For $0 \leq s_1 <... < s_n \leq 1, s=(s_1,...,s_n)$, define $F_s$ to be the convex hull of
$B(s_1),...,B(s_n)$. This is a.s an $n-1$ dimensional simplex. Let $E_s$ be the measure zero event that $F_s$ is a facet in the boundary of the convex hull of the random walk. Our aim is to show that with high probability, none of the events $E_s$ hold for $s_1=0$, which means that the convex hull does not contain any facet  the origin is a vertex of which.  \\

For a point $s$ defined as above, we define $r(s)=(r_1,...r_{n})$ by $r_1 = s_1$, $r_i = s_i - s_{i-1}$ for $2 \leq i \leq n$. The point $r(s)$ lives in the $n$-dimensional simplex, which we denote
by $\Delta_{n}$. Analogously, for a point $r \in \Delta_n$ define by $s(r)$ the corresponding point $s=(s_1,..,s_n)$. By slight abuse of notation we will also write $E_r$ and $F_r$, allowing ourselves to interchange freely between $s$ and $r$. \\

Denote by $W_r$ the measure zero event that the point $r \in \Delta_n$ is also in the poisson process 
(hence the event that all the points $r_1,r_1+r_2,...,r_1+...+r_n$ are in the set $\{0, t_1,...,t_N\}$).  \\ \\
For a Borel subset $A \subset \Delta_{n}$, define
$$
\mu (A) = \EE \left [\sum_{r \in A} \textbf{1}_{E_r} \right ],
$$
the expected number of facets $F_r$, with $r \in A$, and
$$
\nu(A) = \EE \left [\sum_{r \in A} \textbf{1}_{W_r} \right ].
$$
Clearly $\mu$ and $\nu$ are $\sigma$-additive, and $\mu \ll \nu$. Denote 
$$p_n(r) = \frac{d \mu}{d \nu}(r), ~~ \forall r \in \Delta_n.$$
So $p_n(r)$ can be understood as $\P(E_r ~|~ W_r)$. \\ \\
Define $\tilde \Delta_n  = \Delta_n \cap \{r_1=0\}$ and,
$$D = \{ r=(r_1,...,r_n) \in \Delta_n ~~ | ~~ r_i > 0, ~~ \forall 2 \leq i \leq n \}.$$
Let $s = (s_1,...,s_n)$ and $\epsilon > 0$ be such that $s_i - s_{i-1} > \epsilon$ for all $2 \leq i \leq n$. 
Define 
$$
Q = r(\{(x_1,...,x_n); ~ x_i \in [s_i, s_i + \epsilon], ~\mbox{for }i=1,..,n  \}).
$$
Then, by the independence of the number of poisson points on disjoint intervals,
$$
\nu (Q) = \EE \left [\prod_{i=1}^n \#\{j;~t_j \in [s_i, s_i + \epsilon] \} \right ] = (\epsilon \alpha)^n.
$$
By the $\sigma$-additivity of $\nu$, it follows that for a measurable $A \subset \Delta_n \setminus \tilde \Delta_n$,
$$
\nu(A \cap D) = \alpha^n Vol_n(s(A)) = \alpha^n Vol_n(A).
$$
where in the last equality we use the fact that the Jacobian of the function $r \to s(r)$ is identically one. Using analogous considerations on $\tilde \Delta_n$, we get,
$$
\nu(A \cap D) = \alpha^n Vol_n(A) + \alpha^{n-1} Vol_{n-1} (A \cap \tilde \Delta_n)
$$
for all $A \subset \Delta_n$ measurable. By the definition of $p_n(r)$, 
$$
\mu(A) = \alpha^n \int_A p_n(r) d\lambda_n(r) + \alpha^{n-1} \int_{A \cap \tilde \Delta_n} p_n(r) d \lambda_{n-1}(r),
$$
for all measurable $A \subset \Delta_n$, $\lambda_n, \lambda_{n-1}$ being the respective Lebesgue measures. \\ \\

We would like to obtain an upper bound for $\mu (\tilde \Delta_n)$. Using the above formula, this is reduced to obtaining an upper bound for $p_n(r)$. To that end, we use the following idea: the representation theorem for the brownian bridge suggests that we may equivalently construct $B(t)$ by first generating the differences $B(s_j)-B(s_{j-1})$ as independent gaussian random vectors, and then "fill in" the gaps between them by generating a brownian motion up to $B(s_1)$, a brownian bridge for each $1 < j\leq n$, and a "final" brownian motion between $B(s_n)$ and $B(1)$, all of the above independent from each other. To make it formal, fix $r \in \Delta_n$ and define $s = s(r)$. For all $i$, $2 \leq i \leq n$, we write,
$$
D_i = B(t_i) - B(t_{i-1})
$$
and define $C_i:[s_{i-1}, s_i] \to \RR^n$ by,
$$
C_i(t) = B(t) - B(s_{i-1}) - \frac{t - s_{i-1}}{s_i - s_{i-1}} (B(s_i) - B(s_{i-1})),
$$
the bridges that correspond to the intervals $[s_{i-1}, s_i]$. Finally, we define two functions $B_0:[0,s_1] \to \RR^n$ and
$B_f:[s_n,1] \to \RR^n$ by $B_0(t) = B(s_1 - t) - B(s_1)$ and $B_f(t) = B(t) - B(s_n)$. By the independence of the differences of a brownian motion on
disjoint intervals and by the representation theorem for the brownian bridge, it follows that the variables $\{D_i\}_{i=2}^n, \{C_i\}_{i=2}^n, B_0, B_f$ are all
independent, each $C_i$ being a brownian bridge and $B_0$ and $B_f$ being brownian motions. \\ \\
Define $\theta_s$ to be an orthogonal unit normal to $F_s$. Denote,
$$
\tilde C_i = \langle C_i, \theta_s \rangle, ~~\forall 2 \leq i \leq n,
$$
and also $\tilde B_0 = \langle B_0, \theta_s \rangle$ and $\tilde B_f = \langle B_f, \theta_s \rangle$. Since $\theta_s$ is fully determined
by $\{D_i\}_{i=2}^n$, it follows that $\{\tilde C_i\}_{i=2}^n$, $\tilde B_0$ and $\tilde B_f$ are independent. Observe that for all $2 \leq i \leq n$, $\tilde C_i$ is a one-dimensional brownian bridge fixed to be zero at its endpoints, and $B_0$ and $B_f$ are one dimensional brownian motions starting from the origin. \\ \\
A moment of reflection reveals that the event $E_s$ is reduced to the intersection of the following conditions for one of the two possible choices of $\theta_s$: \\ \\
(i) $W_s$ holds. \\
(ii) For all $2 \leq i \leq n$, the function $\tilde C_i$ is non-negative at all points $t_j$ such that $s_i \leq t_j \leq s_{i+1}$. \\
(iii) The function $\tilde B_0$ is non-negative at all points $t_j$ such that $t_j < s_1$. \\
(iv) The function $\tilde B_f$ is non-negative at all points $t_j$ such that $s_n < t_j \leq 1$. \\ \\
As explained above, $\{\tilde C_i\}_{i=2}^n$, $\tilde B_0$ and $\tilde B_f$ are independent, thus we can estimate $p(r)$ using equations (\ref{walkest}) and (\ref{bridgeest}). We get,
\begin{equation} \label{facets}
p_n(r) = \left (\prod_{j=2}^{n} \frac{1}{\alpha r_j} \right ) \frac{1}{\pi} \frac{1}{\sqrt{\alpha r_1} \sqrt{\alpha r_{n+1}}}
\prod_{j=1}^{n+1} \left (1 + O \left (\frac{1}{\alpha r_j} \right ) \right ).
\end{equation}
Using the fact that each probability in the product can be bounded by $1$, we see that there exists a constant $c>0$ such that,
$$
p_n(r) < c^n \left (\prod_{j=2}^{n} \min \left \{ \frac{1}{\alpha r_j}, 1 \right \} \right ) \min \left \{ \frac{1}{\sqrt{\alpha r_1}}, 1 \right \} \min \left \{ \frac{1}{\sqrt{\alpha r_{n+1}}}, 1 \right \} = 
$$
$$
\frac{c^n}{\alpha^n} \left (\prod_{j=2}^{n} \min \left \{ \frac{1}{r_j}, \alpha \right \} \right ) \min \left \{ \frac{1}{\sqrt{r_1}}, \sqrt \alpha \right \} \min \left \{ \frac{1}{\sqrt{r_{n+1}}}, \sqrt \alpha \right \}.
$$
Now,
$$
F(\tilde \Delta_n) = \alpha^{n-1} \int_{\tilde \Delta_n} p(r) d \lambda_{n-1} (r) =
$$
$$
\alpha^{n-1} \int_{\Delta_{n-1}} p_{n-1} (r) \lambda_{n-1} (r) < \alpha^{n-1} \int_{K_{n-1}} p_{n-1}(r) \lambda_{n-1} (r),
$$
where $K_{n-1} = \{0\} \times [0,1]^{n-1}$ is the $n-1$-dimensional cube. So,
$$
F(\tilde \Delta_n) < \alpha^{n-1} \frac{c^n}{\alpha^{n-\frac{1}{2}}} \left ( \int_0^1 \min \{ \frac{1}{r}, \alpha \} dr \right )^{n-1} 
\int_0^1 \min \{ \frac{1}{\sqrt r}, \sqrt \alpha\} dr < 
$$
$$
\frac{c^n}{\sqrt \alpha} \left ( \int_0^1 \min \{ \frac{1}{r}, \alpha\} dr \right )^{n-1} 
\int_0^1 \frac{1}{\sqrt r} dr < \frac{(c' \log \alpha)^n}{\sqrt \alpha}.
$$
Suppose $\alpha = n^{2 L n}$ having $L>3$, then
$$
\frac{(c' \log \alpha)^n}{\sqrt \alpha} = \frac{(2 n L c' \log n)^n}{n^{L n}} =
\left (\frac{2 n L c' \log n}{n^{L}} \right )^n < \left (\frac{2 L c'' }{n^{L - 2}} \right )^n.
$$
We may clearly assume that $n \geq 2$. It follows that there exists a universal constant $C>0$ such that whenever $L \geq C/2$, we have $F(\tilde \Delta_n)<n^{-n}$. Note that the assumption that $L \geq C/2$ may be written $\alpha \geq e^{C n \log n}$. Finally, an application of Markov's inequality then teaches us that in this case, the probability of having one face containing the origin is smaller than $n^{-n}$, which finishes the proof. \qed \\ \\
We have now established theorem \ref{mainttt}. \\
\section{The Discrete Setting} \label{secoedisc}
The aim of this section is to sketch the proof of theorem \ref{discrete}. \\ \\
Fix a dimension $n \in \mathbb{N}$. Let $S_1,...,S_N$ be a standard random walk on $\mathbb{Z}^n$. The following lemma is the discrete analogue of formulas (\ref{walkest}) and (\ref{bridgeest}) derived in the previous section:
\begin{lemma} \label{DiscreteWalk}
Suppose $N>2$. Let $\theta \in \Sph$. Define, 
$$\tilde S_j := \langle \theta, S_j \rangle, ~~ \forall 1 \leq j \leq N.$$
The following estimates hold:
\begin{equation} \label{DiscWalk}
\P \left (\tilde S_j \geq 0, ~~ \forall 1 \leq j \leq N \right ) < \frac{10 n}{\sqrt{N}}
\end{equation}
and,
\begin{equation} \label{DiscBridge}
\P \left (\left . \tilde S_j \geq 0, ~~ \forall 1 \leq j \leq N ~~ \right | ~~ \tilde S_N = 0 \right ) < \frac{2 \log N}{N}.
\end{equation}
\end{lemma}
\textbf{Proof:}
The proof of (\ref{DiscBridge}) follows again from the cyclic shifting principle, explained in the last section. However, it is a bit more involved than the continuous case, since a discrete random walk can attain its global minimum more than once.
Define by $Z_i$ the event that $\tilde S_k = 0$ for exactly $i$ distinct values of $k$, and define,
$$
p_i = \P \left (\left . \left \{ \tilde S_j \geq 0, ~~ \forall 1 \leq j \leq N \right \}~ \cap ~  Z_i ~~ \right | ~~ \tilde S_N = 0 \right )
$$
and,
$$
p = \P \left (\left . \tilde S_j \geq 0, ~~ \forall 1 \leq j \leq N ~~ \right | ~~ \tilde S_N = 0 \right ) = \sum_{i=1}^\infty p_i.
$$
we now use the following observation: consider random walk conditioned on attaining a certain value $T \in \RR$, $\ell$ times. The probability that $T$ is the global minimum of this random walk is smaller than $2^{- \ell}$, since each of the segments between two points can be reflected around the value $T$.  It follows that,
$$
\sum_{i= \lceil \log_2 N \rceil +2}^\infty p_i  \leq \sum_{i=\lceil \log_2 N \rceil + 2 }^\infty 2^{-i+1} \leq \frac{1}{N}.
$$
By the cyclic shifting principle, described in the previous section, we have $p_i \leq i / N$. So,
$$
p = \sum_{i=1}^\infty p_i \leq \frac{1}{N} + \sum_{i=1}^{\lceil \log_2 N \rceil + 2} \frac{i}{N}.
$$
Equation (\ref{DiscBridge}) follows. \\

We turn to prove (\ref{DiscWalk}). Denote $\theta = (\theta_1,...,\theta_n)$. Without loss of generality, we can assume that the $\theta_i$'s are all non-negative and decreasing.
Define the event,
$$
A := \{\tilde S_1 = \theta_1\}.
$$
Clearly,
$$
\P (\tilde S_j \geq 0, ~~ \forall 1 \leq j \leq N) \leq \P(\tilde S_j \geq 0, ~~ \forall 1 \leq j \leq N ~|~ A).
$$
Define $\tilde M_N = \max_{1 \leq j \leq N} \{\tilde S_j\}$. From the symmetry of the random walk,
$$
\P(\tilde S_j \geq 0, ~~ \forall 1 \leq j \leq N ~| ~A) = P(M_{N-1} \leq \theta_1)
$$
Observe that once a random walk went past $\theta_1$ for the first time, it is still at most
$2 \theta_1$. Thus, using the reflection principle, conditioning on the event $M_{N-1} > \theta_1$, we have,
$$
\P(\tilde S_{N-1} > 2 \theta_1 ~|~ M_{N-1} > \theta_1) \leq \frac{1}{2}.
$$
Therefore,
$$
\P(M_{N-1} > \theta_1) \geq 2 \P(\tilde S_{N-1} > 2 \theta_1),
$$
and so,
$$
\P(M_{N-1} \leq \theta_1) \leq 1 - 2 \P(\tilde S_{N-1} > 2 \theta_1) = \P(|\tilde S_{N-1}| \leq 2 \theta_1).
$$
Define,
$$
\phi = (\theta_1,0,...,0) \in \mathbb{R}^n
$$
and define a new random walk, $W_j = \langle \phi, S_j \rangle$. Next we show that for all $a \in \RR$,
\begin{equation} \label{comparewalks2}
\P(|\tilde S_{N-1}| < a) \leq \P(|W_{N-1}| < a).
\end{equation}
Indeed, for all $\lambda \in \RR$,
$$
\EE[ \exp (\lambda \tilde S_{N-1})] = \prod_{j=1}^{N-1} \EE[ \exp (\lambda (\tilde S_{j}-\tilde S_{j-1}))] \geq
$$
$$
\prod_{j=1}^{N-1} \EE[ \exp (\lambda (W_{j}- W_{j-1}))] = \EE[ \exp (\lambda W_{N-1})]
$$
where the last equality follows from the independence of the differences $W_{j}- W_{j-1}$.
Using the symmetry of this differences gives, for all $\lambda \in \RR$,
$$
\EE[ \exp (\lambda \tilde S_{N-1}) + \exp (-\lambda \tilde S_{N-1})] \geq \EE[ \exp (\lambda W_{N-1}) + \exp (-\lambda W_{N-1})],
$$
which implies (\ref{comparewalks2}). We are left with estimating $\P(|W_{N-1}| \leq 2 \theta_1)$. We have,
$$
\P(|W_{N-1}| < a) = \sum_{k=0}^{N-1} \left ( \frac{1}{n}^k \left (\frac{n-1}{n} \right )^{N-1-k} \left (N-1 \atop k \right)\sum_{j=-2}^2 \left ( k \atop \lfloor \frac{k}{2} \rfloor + j \right ) \right ) < \frac{10 n}{\sqrt N}.
$$
This finishes the proof. \qed \\ \\
\textbf{Sketch of the proof of theorem \ref{discrete}}:
We begin with the upper bound. We follow that same lines as the ones in the proof of theorem \ref{upper}. The only extra tool needed for the proof of the upper bound is lemma \ref{DiscreteWalk}. \\ \\ Fix $N \in \mathbb{N}$. For $1 \leq j \leq N$ and $t=\frac{j}{N}$, define $B(t) := Z_{j}$. Let $r=(r_1,...,r_n) \in \Delta_n \cap \frac{1}{N} \mathbb{Z}^n$, and $t_k = \sum_{j=1}^k r_j$. Define the event $E_r$ in the same manner: 
$$
E_r := \{ conv(B(t_1),...,B(t_n)) \mbox{ is contained in the boundary of } K \}
$$
For $A \subset \Delta_n \cap \frac{1}{N} \mathbb{Z}^n$, define
$$F(A) = \EE \left [ \sum_{r \in A} \mathbf{1}_{\{ E_r \}} \right ].$$
Next, for any $r \in \Delta_n \cap \frac{1}{N} \mathbb{Z}^n$, equations (\ref{DiscWalk}) and (\ref{DiscBridge}) are used to obtain,
$$
\P(E_r) < 100 (\log N)^{2n} n^2 \left (\prod_{j=2}^{n} \min \left \{ \frac{1}{N r_j}, 1 \right \} \right ) \min \left \{ \frac{1}{\sqrt{Nr_1}}, 1 \right \} \min \left \{ \frac{1}{\sqrt{N r_{n+1}}}, 1  \right \}.
$$
Define $\Delta_0 = \Delta_n \cap \frac{1}{N} \mathbb{Z}^n \cap \{r_1=0\}$. 
We are left with estimating,
$$ F(\Delta_0) = \sum_{r \in \Delta_0} \P(E_r).$$
This can be done by showing that these are Riemann sums converging to an integral which
can be estimated in the same manner as in theorem \ref{upper}. An analogous calculation gives,
$$
F(\Delta_0) \leq \frac{(C n^2 \log^3 N)^n}{\sqrt{N}}.
$$
For some universal constant $C>0$, which implies the upper bound.
$$~$$
Next, we prove the lower bound. Again follow the same lines as in the proof of theorem \ref{lower}. \\ \\
Assume that $N= 2^{m-1}$ where $m = \left \lfloor \frac{c n}{\log n} \right \rfloor$, the value of the 
constant $c$ will be chosen later. We construct a vector $v$ in an analogous manner to the construction in theorem \ref{lower}. Define $v_0 = Z_1$ and,
$$
v_i = S_{2^i} - S_{2^{i-1}} 
$$
for $i=1,...,m-1$. Define,
$$
v = \frac{1}{\sqrt m} \sum_{i=0}^{m-1} \frac{v_i}{\sqrt{\EE[|v_i|^2]}} =  \frac{1}{\sqrt m} \sum_{i=1}^m \frac{v_i}{(\sqrt{2})^{i-1}}
$$
Fix a $1 \leq k \leq m$, and define 
$$
p = \langle S_{2^k}, \theta \rangle.
$$
The expectation and variance of $p$ can be computed directly, as in the proof
of theorem \ref{lower}. Defining, the $w_j$'s analogously, Chernoff's inequality can be used to prove the bound (\ref{moments}). Theorem \ref{berns} is used to show that for a small enough value of $c$,
$$
\P(p < \frac{1}{2} \EE[p]) < n^{-5}.
$$
By applying a union bound, we can make sure that $\langle S_{2^k}, \theta \rangle \geq \frac{1}{2} \EE \left [ \langle S_{2^k}, \theta \rangle \right ]$ for all $1 \leq k \leq m$.  Next, a formula analogous to (\ref{bridgeprop}) should be applied in order to control the conditional random walks found between consecutive points of the form $2^k$. To this end, we observe that for our random walk $\tilde S_n := \langle \theta, S_n \rangle$ one has,
$$
\PP \left ( \max_{1 \leq j \leq k} \tilde S_j < u \right ) \leq \PP \left ( \left. \max_{1 \leq j \leq k} \tilde S_j < u ~ \right |  \tilde S_k = 0 \right ), ~~ \forall k \in \mathbb{N}, u>0.
$$
Hence, instead of bounding a conditional random walk, we may bound the usual random walk. Using Bernstein's 
inequality, theorem \ref{berns}, in order to derive a bound analogous to (\ref{bridgebound}). Using a union bound gives,
$$
\P(\langle S_j, v \rangle > 0, \forall 1 \leq j \leq N) > 1 - \frac{1}{n}.
$$
This finishes the sketch of proof. \qed 

\section{Spherical covering times} \label{covsec}
The goal of this section is to prove corollary \ref{covering}. \\ \\

Let $B(t)$ be a standard brownian motion in $\RR^n$, $n>2$. Denote $\theta(t) = \frac{B(t)}{|B(t)|}$ and observe that $\theta(t)$ is almost surely
well-defined for all $t>0$. Let $T(t)$ be the solution of the equation 
$$
T'(t) = |B(T(t))|^2, ~~ T(0)=1.
$$
We denote by $[S]_t$ the quadratic variation of an It\"o process, $S_t$, between time $0$ and time $t$. We have,
$$
\frac{d}{dt} [\theta \circ T]_t = T'(t) \frac{(n-1) \left ( \frac{d}{dt} [B]_{t} \right ) |_{t=T} }{n |B(T(t))|^2} = 
$$
$$
\left . \frac{n-1}{n} \left ( \frac{d}{dt} [B]_{t} \right ) \right |_{t=T} = n-1.
$$
which implies that $\theta(T(t))$ is a strong Markov process, and is therefore a spherical brownian motion. \\ \\

\emph{Proof of corollary \ref{covering}}: \\
First, observe that for every $\tau > 0$, the origin lies in the interior of $conv(\{ B(t) ; 1 \leq t \leq \tau \})$ if and only if it lies in the interior of $conv(\{ \theta(t) ; 1 \leq t \leq \tau \})$, thus we have $E(n) = \EE [ \tau_1 ]$ where
$$
\tau_1 = \inf \left \{\tau > 0; ~ F_\tau \mbox{ holds} \right \},
$$
and
$$
F_\tau = \{0 \in Int(conv(\{B(T(s)) ; ~ 0 \leq s \leq \tau\})) \}.
$$
$$~$$
We aim to use the bounds from theorems \ref{lower} and \ref{upper}. For that,
we will need to establish certain bounds on the distribution of $T^{-1}(s)$ for a given $s>0$. \\ \\
Since $\EE(|B(T)|^2) = nT$, it follows that $\EE(T(t)) = e^{nt} + 1$. Using Markov's inequality gives
\begin{equation} \label{eq848837}
\PP \left (T(t) > 10 e^{nt} + 10 \right ) \leq 0.1.
\end{equation}
By theorem \ref{lower}, there exists a constant $c>0$ such that for
$$
\tau_2 = \inf \{ \tau > 0; ~ T(\tau) \geq e^{c n / \log n}  \},
$$
one has
\begin{equation} \label{ftau2small}
\PP(F_{\tau_2}) < 0.1.
\end{equation}
According to equation (\ref{eq848837}), 
\begin{equation} \label{tau2small}
\PP(\tau_2 < c_1 / \log n) < 0.1,
\end{equation}
for some universal constant $c_1 > 0$. Using a union bound with (\ref{ftau2small}) and (\ref{tau2small}) gives,
$$
\PP(\tau_1 <  c_1 / \log n) < 0.2,
$$
which implies
$$
\EE[\tau_1] \geq 0.8 c_1 / \log n.
$$
The lower bound is established. \\ \\
We continue with the upper bound. Observe that $T(t)$ is a bijective map from $[0,\infty)$ to $[1,\infty)$. We may define $f(s) = T^{-1}(s)$ for all $s \geq 1$. One has,
$$
f'(s) = \frac{1}{T'(f(s))} = \frac{1}{|B(s)|^2}.
$$
Consequently, by Fubini's theorem,
$$
\EE[f(s)] = \int_1^s \EE \left [ \frac{1}{|B(t)|^2} \right ] dt = \int_1^s \frac{1}{t} \EE \left [ \frac{1}{|\Gamma|^2} \right ] dt,
$$
where $\Gamma$ is a standard gaussian random vector in $\RR^n$. A calculation gives $\EE \left [ \frac{1}{|\Gamma|^2} \right ] < \frac{C_1}{n}$
for some universal constant $C_1>0$. It follows that $\EE[f(s)] \leq \frac{C_1 \log s}{n}$. By Markov's inequality,
\begin{equation} \label{Ttlowerbound}
\PP \left (f(s) > \frac{10 C_1 \log s}{n} \right ) < 0.1.
\end{equation}
According to theorem \ref{upper}, there exists a universal constant $C>0$ such that for
$$
\tau_3 = \inf \{ \tau > 0; ~ T(\tau) \geq e^{C n \log n}  \},
$$
one has,
\begin{equation} \label{ftau2small2}
\PP(F_{\tau_3}) > 0.9.
\end{equation}
Now, an application of equation (\ref{Ttlowerbound}) with $s = e^{C n \log n}$ gives,
\begin{equation} \label{ftau2small3}
\PP(\tau_3 > C_2 \log n) < 0.1.
\end{equation}
for some universal constant $C_2>0$. Using a union bound with equations (\ref{ftau2small2}) and (\ref{ftau2small3}) gives,
$$
\PP(\tau_1 >  C_2 \log n) < 0.2.
$$
In other words,
$$
\PP \left ( 0 \in Int(conv(\{ \theta(T(t)) ; 0 \leq t \leq C_2 \log n \})) \right ) > 0.8.
$$
Now, by the strong Markov property and time-homogeneouity of $\theta \circ T$, we also have
$$
\PP \left ( 0 \in Int(conv(\{ \theta(T(t)) ; k C_2 \log n \leq t \leq (k+1) C_2 \log n \})) \right ) > 0.8.
$$
for all $k \in \mathbb{N}$. Finally, since the above event is invariant under rotations,
$$
\PP \left ( 0 \in Int(conv(\{ \theta(T(t)) ; 0 \leq t \leq k C_2 \log n \})) \right ) > 1 - 0.2^k.
$$
In other words,
$$
\PP(\tau_1 > C_2 k \log n) < 0.2^k,
$$
which easily implies that $E[n] \leq C_3 \log n$, for some universal constant $C_3 > 0$. 
The proof is complete. \qed \\ \\

\section{Remarks and Further Questions} \label{secoeremarks}
In this section state a few results that can easily be obtained using the same ideas used above, and suggest possible related directions of research.
\subsection{Probability for intermediate points in the walk to be extremal.}
The methods used above can easily be adopted in order to estimate the probability that an intermediate point of a random walk is an extremal point. To see this, observe that this probability is equivalent to the probability that the origin is an extremal point of two independent random walks of length $\lambda N$ and $(1-\lambda)N$ respectively. Thus, theorem \ref{upper} can still be used for an upper bound since either $\lambda \geq \frac{1}{2}$ or $1-\lambda \geq \frac{1}{2}$. For the lower bound we should do a little extra work: we follow the lines of the proof of theorem \ref{lower}, only defining the vector $v$ as,
$$
v = \lambda v_1 + (1 - \lambda) v_2
$$ 
where $v_1$ and $v_2$ are constructed in the same manner that the vector $v$ is constructed in theorem \ref{lower}. The exact same calculations can be carried out to show that with high probability $v$ separates the origin from the points of both of the random walks. This yields,
\begin{proposition}
There exist universal constants $C,c>0$ such that the following holds:
Let $S_1,S_2,...$ be the standard random walk on $\mathbb{Z}^n$ and let $j, N \in \mathbb{N}, ~~ j < N$. Then: \\
(i) If $N > e^{C n \log n}$ then $\PP( S_j \in Int(conv \{S_1,...,S_N\})) > \frac 1 2$. \\
(ii) If $N < e^{c n / \log n}$ then $\PP( S_j \in \partial conv \{S_1,...,S_N\}) > \frac 1 2$.
\end{proposition}

\subsection{Covering times and Comparison to independent origin-symmetric random points}
The result of corollary \ref{covering} can also be viewed as an upper bound on a certain mixing time of the spherical brownian motion: 
Let $\mu$ be an origin-symmetric distribution on $\RR^n$ which is absolutely continuous with respect to the Lebesgue measure. There is a beautiful proof by Wendel, \cite{W}, if $X_1,...,X_N$ are independent random vectors with law $\mu$, one has
\begin{equation} \label{indpoints}
\P(0\notin \mbox{conv}\{X_1,\dots,X_N\})=\frac{1}{2^{N-1}}\sum\limits_{k=0}^{n-1}{N-1 \choose k }. 
\end{equation}
Hence, the probability does not depend on $\mu$ as long as it is centrally symmetric and absolutely continuous. Note that in order for this probability to be $\frac{1}{2}$ one should take $N(n) \approx n \log n$. \\
This suggests that the correct mixing time in the sense of the $\frac{\pi}{2}$-covering should be $\frac{1}{n}$. \\ \\

An easy computation shows that after time of order $\frac{1}{n}$, a brownian motion that started at an arbitrary point on 
the sphere will be approximately uniformly distributed on the sphere in the sense that the density will be bounded between two universal constants, independent of the dimension. If we assume that the correct mixing time is therefore $\frac{1}{n}$
for this purpose, this suggests that our upper bound of $e^{n \log n}$ should be a natural conjecture for the correct asymptotics in theorem \ref{mainttt}.

\subsection{A random walk that does not start from the origin}
Our techniques may be also used to find the asymptotics of the time it takes for the origin to be encompassed by a random walk when the 
starting point is different than the origin. By the scaling property of brownian motion,
$$
\PP(0 \in Int(Conv\{B(t) ; 1 \leq t \leq M \})) = \PP(0 \in Int(Conv\{B(t) ; L \leq t \leq L M \}))
$$
For all $M>1, L>0$. Using the concentration of $|B(t)|$ around its expectation, it is not hard to derive,

\begin{proposition}
There exist universal constants $C,c>0$ such that the following holds:
Let $B(t)$ be a brownian motion started at a point $x_0$ whose distance from the origin is $L$. Then: \\
(i) If $M > L^2 e^{C n \log n}$ then $\PP(0 \in Int(conv \{B(t); ~ 0 \leq t \leq M \})) > \frac 1 2$. \\
(ii) If $M < L^2 e^{c n / \log n}$ then $\PP(0 \in Int(conv \{B(t); ~ 0 \leq t \leq M \})) < \frac 1 2$.
\end{proposition}

\subsection{Possible Further Research}

In this note we try to find the correct asymptotics, with respect to the dimension $n$, of the value $N$ such that $p(n,N) \approx \frac{1}{2}$.
One related question is:
\begin{question}
For a fixed value of $n$, how does $p(n,N)$ behave asymptotically as $N \to \infty$?
\end{question}

In view of (\ref{indpoints}) and the discussion following it, 
one might expect that this probability could have approximately the following law, 
for a certain range of values of $N$,
$$
p \approx \frac{(\log N)^n}{N^c}
$$
where $p$ is the probability in question, $n$ is the dimension and $N$ is the length
of the random walk, and $c>0$ is some constant. \\ \\
Two other possible questions are:
\begin{question}
Given two numbers $j,k < N$, what is the joint distribution of $S_j, S_k$ being extremal points
of the random walk $S_1,...,S_N$? Is there repulsion or attraction between extremal points of
a random walk?
\end{question}

\begin{question}
How does the result of theorem \ref{mainttt} change is one replaces the brownian motion by a $p$-stable process?
\end{question}

\bigskip {\noindent School of Mathematical Sciences, Tel-Aviv University, Tel-Aviv
69978, Israel \\  {\it e-mail address:}
\verb"roneneldan@gmail.com" }

\end{document}